\documentclass[12pt]{amsart}
\usepackage{amscd,amssymb}
\usepackage{graphicx,enumerate}

\usepackage{color}
\usepackage[arrow,matrix,arrow,ps,color,line,curve,frame]{xy}

\usepackage{pgf, tikz}
\usepackage{url}
\usepackage{enumerate}

\usepackage[colorinlistoftodos,bordercolor=red,backgroundcolor=red!20,linecolor=red,textsize=scriptsize]{todonotes}

\def\height{\operatorname{ht}} 

\def\width{\operatorname{width}} 
\def\conv{\operatorname{conv}}

\def\inte{\operatorname{int}}
\def\Hilb{\operatorname{Hilb}}
\def\Pl(d){\operatorname{Pol}(d)} 

\def\vol{\operatorname{vol}}

\def\np{{\operatorname{NPol}}}

\def\k{{\mathbf k}}

\def\POL{{\mathbf{POL}}}

\def\cones{{\textsf{Cones}}} 
\def\RR{{\mathbb R}} 
\def\QQ{{\mathbb Q}}
 
\def\ZZ{{\mathbb Z}} 
\def\NN{{\mathbb N}}
\def\icp{\operatorname{ICP}} 
\def\ucp{\operatorname{UCP}} 
\def\Aff{\operatorname{Aff}}

\let\Bbb=\mathbb

\let\epsilon=\varepsilon

\def\CR{\operatorname{CR}}

\def\cones{\operatorname{\textsf{Cones}}}

\def\width{\operatorname{width}}

\def\vol{\operatorname{vol}}%

\def\Hilb{\operatorname{Hilb}}%

\def\Q{{\Box\kern1pt}}
\def\inte{\operatorname{int}}

\def\D{\underset\Delta}
\def\q{\underset{\bf q}}
\def\Infty{\underset{\infty}}

\def\conv{\operatorname{conv}}%
\def\vol{\operatorname{vol}}%

\def\RR{{\Bbb R}}

\def\QQ{{\Bbb Q}}
\def\ZZ{{\Bbb Z}}
\def\NN{{\Bbb N}}
\def\ee{{\mathbf e}}
\def\zz{{\mathbf z}}

\newtheorem{lemma}{Lemma}[section]  \newtheorem{theorem}[lemma]{Theorem}
 \newtheorem{conjecture}[lemma]{Conjecture} 

\theoremstyle{definition} \newtheorem{definition}[lemma]{Definition}

\advance\headheight1.15pt \marginparwidth=2.5cm

\begin{document}

\title[Normal polytopes]{Normal polytopes: between discrete,\\ continuous, and random}

\author[Joseph Gubeladze]{Joseph Gubeladze}

\address{Department of Mathematics, San Francisco State University, San Francisco, CA 94132, USA}

\email{soso@sfsu.edu}


\begin{abstract}
The first three sections of this survey represent an updated and much expanded version of the abstract of my talk at FPSAC'2010: new results are incorporated and several concrete conjectures on the interactions between the three perspectives on normal polytopes in the title are proposed. The last section outlines new challenges in general convex polytopes, motivated by the study of normal polytopes.
\end{abstract}

\maketitle

\section{Normal polytopes: general overview}\label{General}

All our polytopes are assumed to be convex.

The group of affine automorphisms of $\RR^d$, mapping $\ZZ^d$ bijectively to itself, will be denoted by $\Aff(\ZZ^d)$.

Let $P\subset\RR^d$ be a \emph{lattice polytope,} i.e., with vertices in $\ZZ^d$. Denote by $\Lambda$ the
subgroup of $\ZZ^d$, generated by the lattice points in $P$.

Here is the central definition:

(a) $P$  is \emph{integrally closed} if the following condition is
satisfied:
$$
c\in\NN,\ z\in(cP)\cap\ZZ^d\quad\Longrightarrow\quad\exists
x_1,\ldots,x_c\in P\cap\ZZ^d\quad x_1+\cdots+x_c=z.
$$

(b) $P$ is \emph{normal} if the following condition is satisfied
$$
c\in\NN,\ z\in(cP)\cap\Lambda\quad\Longrightarrow\quad\exists
x_1,\ldots,x_c\in P\cap\Lambda\quad x_1+\cdots+x_c=z.
$$

The normality property is invariant under the affine
isomorphisms of lattice polytopes, respecting the sets of lattice points in the polytopes, and `integrally closed\rq{} is invariant under the $\Aff(\ZZ^d)$-action.

Obvious examples of normal but not integrally closed polytopes are
the \emph{empty lattice simplices} of large volume. A full classification of such simplices is only known in dimensions $\le4$ \cite{Empty4}.

A normal polytope $P\subset\RR^d$ can be made into an
integrally closed polytope by changing $\RR^d$ to the subspace $\RR\Lambda$ and the lattice of reference $\ZZ^d$ to $\Lambda$. In particular, normal and integrally closed polytopes refer to same isomorphism classes of lattice polytopes. In the combinatorial literature the difference between `normal' and `integrally closed' is sometimes blurred.

\medskip\noindent\emph{Convention.} In this text we decided to follow the more recent practice of using `normal\rq{} in the sense of `integrally closed\rq{}. 

\medskip Normal polytopes form an important  subclass of lattice polytopes, central in combinatorial commutative algebra and toric algebraic geometry: they define the normal homogeneous monoid algebras  \cite{kripo,MiStcomcom} and the projectively normal toric embeddings \cite{kripo,CLSTor}. These polytopes are also important in integer programming \cite{Schrijver}. Normal polytopes enjoy popularity in algebraic combinatorics and they have been showcased on several workshops \cite{aim,mfo}. There are many challenges of number theoretic, ring theoretic, homological, and $K$-theoretic
nature, concerning the associated objects: rational cones, affine monoid rings, and toric varieties \cite{kripo}.

Our approach to normal polytopes in this survey is not exactly along the lines of the mentioned connections, though. Our starting point is that normal polytopes are natural discrete models of convex compact sets: informally, the normality condition says that the distribution of lattice points stays homogeneous no matter how far we zoom into the polytope. Correspondingly, the goal is to understand how far this continuous/discrete analogy can be promoted. Of course, the idea that general lattice polytopes are discrete models of general convex polytopes is not new (e.g., \cite{Beck}), but our emphasize is exclusively on normal polytopes: the questions we raise below trivialize without the extra condition of normality.

If a lattice polytope is a union of (in particular, subdivided into) normal polytopes then it is also normal. The simplest normal polytopes one can think of are \emph{unimodular simplices}, i.e., the lattice simplices $\Delta=\conv(x_1,\ldots,x_k)\subset\RR^d$, $\dim\Delta=k-1$, with $x_1-x_j,\ldots,x_{j-1}-x_j,x_{j+1}-x_j,\ldots,x_k-x_j$ a part of a basis of $\ZZ^d$ for some (equivalently, every) $j$. The next simplest class of normal polytopes is given by \emph{lattice parallelepipeds} \cite[\S3]{unico}.

A classical observation is that every lattice polytope $P\subset\RR^d$ of dimension $\le2$ admits a unimodular triangulation and, therefore, is normal \cite[Corollary 2.54]{kripo}. However, starting from dimension $d\ge3$, the majority of lattice polytopes are not normal.

Unimodular simplices are the smallest `atoms' in the world of
normal polytopes. But the latter is not built out exclusively of
these atoms, i.e., not all normal polytopes are unions of unimodular simplices. This and other `negative' results contributed to the
current thinking in the area that there is no succinct geometric
characterization of the normality property. One could even
conjecture that in high dimensions the situation gets as bad as
it can; see Section \ref{CR} for details.

We say that a lattice polytope has the \emph{unimodular cover property} ($\ucp$) if it is a union of unimodular simplices.

It is classically known \cite[Chapter III]{kkms} that every lattice polytope $P$ has a multiple $cP$ for some $c\in\NN$ that is \emph{triangulated} into unimodular simplices. The recent preprint \cite{Gaku} strengthens this result to the existence of a threshold value $c$, depending on $P$, for which \emph{all} multiples $c'P$, $c'\ge c$, have unimodular triangulation. The existence of a \emph{dimensionally uniform} lower bound for such $c$ is a major open problem. 

The situation is better for $(\ucp$). Namely, an explicit dimensionally uniform lower bound for the factors $c$, for which the multiples $c'P$ with $c\rq{}\ge c$ have ($\ucp$), was derived in \cite{multiples}. Later, by improving one important step in \cite{multiples}, the original exponential bound was cut down to a degree 6 polynomial function in the dimension \cite[\S3C]{kripo}, \cite{vonthaden}.

These results on unimodular covers of multiple polytopes yield no new examples of normal polytopes, though. In fact, an easy argument ensures that, for any lattice $d$-polytope $P$, the multiples $cP$ with $c\ge d-1$ are normal \cite[Prop. 1.3.3]{brgutr}, \cite{ewald}. One should remark that there is no algebraic obstruction to the existence of very special -- the so called \emph{regular quadratic} unimodular triangulations for the multiples $cP$ with $c\ge d$: the nice homological
properties that the corresponding toric rings would have
according to \cite{StuTOH} in the presence of such triangulations can be independently derived via alternative algebraic methods \cite{brgutr}.

Lattice polytopes with long edges of independent lengths are
considered in \cite{edges}, where it is shown that if the edges of
a lattice $d$-polytope $P$ have lattice lengths $\ge 4d(d+1)$
then $P$ is normal; moreover, in the special case when $P$ is a
simplex one can do better: $P$ is covered by lattice
parallelepipeds, provided the edges of $P$ have lattice lengths
$\ge d(d+1)$. Later the bound $4d(d+1)$ was improved to $2d(d+1)$ in \cite{halved}.

Currently, one problem attracts much attention in the field.
Namely, \emph{Oda's question} asks whether  all \emph{smooth} polytopes are normal. A lattice polytope $P\subset\RR^d$ is called
smooth if the primitive edge vectors at every vertex of $P$ define
a part of a basis of $\ZZ^d$. Smooth polytopes correspond to
projective embeddings of smooth toric varieties. Oda's
question is open in all dimensions $\ge3$. 

Recent positive results in the field mostly concern special classes of polytopes: associated to graphs or root systems, of low dimension, satisfying additional geometric constrains, or occurring in various applications. But this is outside of the scope of our survey.

\section{Unimodular cover and Carath\'eodory rank}\label{CR}

There are 4-dimensional normal polytopes without unimodular triangulations \cite[Proposition 1.2.4.(c)]{brgutr}. 

At the end of the 1980s, it was conjectured \cite{cfs} that every rational finite pointed cone has the s.c. \emph{integral Carath\'eodory property}. Via the \emph{homogenization} (see Section \ref{Spaces}) normal polytopes can be thought of a special case of rational cones. The mentioned conjecture in the special case of a normal polytope $P$ says that for every $c\in\NN$ and every $x\in(cP)\cap\ZZ^d$, there exist points $x_1,\ldots,x_{d+1}\in P\cap\ZZ^d$ and numbers $a_1,\ldots,a_{d+1}\in\ZZ_{\ge0}$ such that $a_1+\cdots+a_{d+1}=c$ and $a_1x_1+\cdots+a_{d+1} x_{d+1}=x$. Let ($\icp$) denote the mentioned property of lattice polytopes. A stronger conjecture was proposed \cite{SebHILB}, according to which all normal polytopes have ($\ucp$). Obviously, for every lattice polytope, ($\ucp$) impllies ($\icp$). In its turn, it can be shown that ($\icp$) implies normality \cite[Theorem 6.1]{unico}. A five-dimensional example of a normal polytope without ($\ucp$) was found in \cite{unico}. The strategy for generating candidates for a counterexample to ($\ucp$) was to first randomly generate lattice parallelepipeds and then construct descending sequences of normal polytopes of the form discussed in Section \ref{Poset}. This approach is based on the observation \cite[Corollary 2.3]{unico} that the validity of the implications

\medskip\centerline{\fbox{Normality\ \ \ \ $\Longrightarrow$\ \ \ \ ($\ucp$)}\ \ and \ \ \fbox{Normality\ \ $\Longrightarrow$\ \ ($\icp$)}}

\medskip\noindent can be checked on the terminal objects of such descending sequences. Then it was shown in \cite{cara} that the discovered 5-dimensional counterexample to $(\ucp$) also violates ($\icp$). That   $(\ucp)$ and $(\icp)$ are in fact two different properties was shown considerably later in \cite{icp}, utilizing the idea that the stronger property ($\ucp$) is likely to be lost before ($\icp$) along the mentioned descending chains.

After many years of extensive computer search only a handful of counterexamples to ($\ucp$) or ($\icp$) have emerged, the majority of which are closely related to the initially found polytope.

For a lattice polytope $P$, we define its \emph{Carath\'eodory rank} $\CR(P)$ as the smallest natural number $k$ such that for every natural number $c$ and every lattice point $z\in cP$ there exist lattice points $x_1,\ldots,x_k\in P$ and integers
$a_1,\ldots,a_k\ge0$ such that $z=a_1x_1+\cdots+a_kx_k$ and
$a_1+\cdots+a_k=c$.

The Carath\'eodory rank allows to quantify the failure of ($\icp$). For every normal polytope $P$, one has
$$
\dim P+1\le\CR(P)\le 2\dim P,
$$
where the first inequality is straightforward and the second follows from \cite{SebHILB}. For $d\in\NN$, put
\begin{align*}
\CR(d)=\max(\CR(P)\ |\ P\subset\RR^d\ \text{a normal polytope}\},\qquad d\in\NN.
\end{align*}
We conjecture that the upper estimate for $\CR(P)$ above is asymptotically sharp:
\begin{conjecture} $\frac{\CR(d)}d$ is monotonically converging to $2$ as $d\to\infty$.
\end{conjecture}
The counterexample to (ICP) in \cite{cara} yields the estimate
$$
\limsup_{d\to\infty}\frac{\CR(d)}d\ge\frac76.
$$

The conjecture in particular says that there are essentially new counterexamples to ($\icp$) in higher dimensions, not obtained by trivial extensions of the counterexamples in lower dimensions. 

For a normal polytope, its Carath\'eodory rank can be computed algorithmically. This follows from \cite[\S4]{BrGuTrPROB}, where a general algorithm is described, without implementation, for checking whether a finitely generated submonoid $M\subset\ZZ^d$ is the union of given finitely generated submonoids $M_1,\ldots,M_n\subset M$. In particular, for a normal polytope $P$, the inequality $\CR(P)\le r$ can be verified by checking whether the submonoid $\RR_{\ge0}(P,1)\cap\ZZ^{d+1}\subset\ZZ^{d+1}$ is the union of its submonoids of the form $\ZZ_{\ge0}(x_1,1)+\cdots+\ZZ_{\ge0}(x_r,1)$, where $x_1,\ldots,x_r\in P\cap\ZZ^d$.

In lieu of the missing systematic procedure for deriving counterexamples to ($\icp$), we propose to take the tried random search to the extreme: one possible explanation of the scarce counterexamples to ($\icp$) or ($\ucp$) so far might be an inherent statistical bias in the pseudo-random numbers, used in the computer; see Section \ref{Search} for details.

\section{Quantum jumps}\label{QJumps}

\subsection{Poset of normal polytopes}\label{Poset} The set of normal polytopes in $\RR^d$, denoted by $\np(d)$, is partially ordered as follows: $Q<P$ if and only if there is a finite sequence of normal polytopes of the form:
\begin{align*} 
&Q=P_0\subset\ldots\subset
P_{n-1}\subset P_n=P,\\ 
&\#(\ZZ^d\cap P_{i+1})=\#(\ZZ^d\cap P_i)+1,\quad i=0,\ldots,n-1. 
\end{align*}

When $\dim(P_i)=\dim(P_{i-1})+1$ then $P_i$ is a \emph{unimodular pyramid} over $P_{i-1}$.

As described in Section \ref{CR}, minimal elements in $\np(5)$ proved crucial in disproving $(\icp)$.

The poset $\np(d)$ is a discrete version of the continuum of all convex compact sets in $\RR^d$, ordered by inclusion. In this context, it provides a formalism for the  dichotomy `discrete vs. continuous' and allows to study global properties of the family of normal polytopes in analogy with \emph{moduli spaces,} as opposed to studying particular (classes of) polytopes. The distortion of the continuum of convex compacta in $\RR^d$ in its discrete model $\np(d)$ is encoded in the homotopy type of the order simplicial complex of the latter.

We introduce the following notation. Every lattice $d$-polytope $P\subset\RR^d$ gives rise to a stratification of the set $\ZZ^d\setminus P$ as follows. Let $P=\bigcap_{i=1}^n H_i^+$ be the irredundant representation as the intersection of affine half-spaces, whose boundary affine hyperplanes are the $H_i$. The lattices $H_i\cap\ZZ^d\subset H_i$, $i=1,\ldots,n$, are bijective to $\ZZ^{d-1}$ via affine maps. For an index $i\in\{1,\ldots,n\}$ and a natural number $j$, denote by $H_i^+(-j)$ the parallel translate of $H_i^+$, such
that $\ZZ^d\cap\left(H^+_i(-j)\setminus H_i^+\right)$ is the union of $j$ parallel translates of $H_i\cap\ZZ^d$, i.e., $H_i(-j)$ is the parallel translate of $H_i$ on lattice distance $j$ \emph{away from} $P$.  One has the stratification:
\begin{align*}
\ZZ^d\setminus P=\bigcup_{j=1}^\infty\big(\partial(
P^{-j})\cap\ZZ^d),\quad\text{where}\ \ &P^{-j}=\bigcap_{i=1}^nH_i^+(-j),\quad j=1,2,\ldots,\ \ \text{and}\\
&\partial(-)\ \ \text{refers to the boundary}.
\end{align*}

Informally, $\partial(P^{-j})\cap\ZZ^d$ consists of the lattice points `on lattice distance $j$ from $P$\rq{}. The polytopes $P^{-j}$ are rational, usually not lattice. 

A \emph{quantum jump} in $\np(d)$ is a pair $(P,Q)$ of $d$-dimensional polytopes, forming an elementary relation $P<Q$ in $\np(d)$. The adjective `quantum\rq{}, apart of a smallest possible change, points to random walks on $\np(d)$: among all possible quantum jumps one can choose the ones according to an adopted strategy. This can be done by introducing various measures: for a quantum jump $(P,Q)$, one defines its (i) \emph{height} by $\height(P,Q)=j$, where $v$ is the vertex of $Q\setminus P$ and $v\in\partial P^{-j}$, and (ii) \emph{volume} $\vol(P,Q)$ as the volume  of the set $Q\setminus P$, normalized with respect to $\ZZ^d$. The volume and height of a quantum jump are always natural numbers.

Below are some of the properties of quantum jumps from \cite{Jumps}:
\begin{enumerate}[\rm$\centerdot$]
\item For any quantum jump $(P,Q)$ in $\np(d)$, one has $\height(P,Q)\le1+(d-2)\width(P)$, where `$\width$\rq{} stands for the the maximum of the lattice widths with respect to the facets; consequently, there are only finitely many quantum jumps from a normal $d$-polytope;
\item For every $d\ge3$, there is no uniform upper bound for the heights of quantum jumps in $\np(d)$;
\item After running many millions of examples, no maximal or $3$-dimensional minimal element has been found in $\np(3)$, although it is known that the order in $\np(3)$ is different from the inclusions order \cite[Example 2.4]{Cones};
\item $\np(4)$ and $\np(5)$ contain maximal elements.
\end{enumerate}

It should be mentioned that, for every $d\ge4$, there are $d$-dimensional minimal elements of $\np(d)$. In fact, for any 4-dimensional minimal element $P\in\np(4)$ and every natural number $d\ge4$, the polytope $P\times[0,1]^{d-4}$ is a minimal element of $\np(d)$. Whether there are maximal elements in $\np(d)$ for all $d$ is not known -- there is no obvious way to derive high-dimensional maximal normal polytopes from lower dimensional ones. 

The above mentioned maximal element in $\np(4)$ and $\np(5)$ were found by the following two independent methods, relevant also in Section \ref{Search} below:

\begin{enumerate}[\rm$\centerdot$]
\item \emph{(Greedy)} Moving along ascending chains of normal polytopes, maximizing at each step the volume of the quantum jump: this approach randomizes well the directions in which the new vertices are added successively; the opposite type of polytope growth in $\np(d)$ occurs when the new vertices align in one direction and that leads to an infinite sequence of similar quantum jumps;
\item \emph{(Random)} Randomly generating lattice polytopes, where the fast check of the normality condition becomes even more decisive -- something \textsf{Normalize} \cite{Normaliz} accomplishes par excellence.
\end{enumerate}

We like to think of $\np(d)$ in terms of the following physical concepts, which also motivates our approach to this global object: the height of a quantum jump represents the \emph{duration} and the volume represents the \emph{energy} of the jump. Informally, the minimal elements of $\np(d)$ model \emph{initial states}, popping into existence, and the maximal elements model \emph{terminal states}. The isolated elements of $\np(d)$, if they exist at all, would  model \emph{eternal states}. The isolated normal polytopes represent the most extremal lattice point configurations, also quite interesting from the additive number theory perspective. The conjecturally nontrivial homotopy type of the order simplicial complex of $\np(d)$ can be interpreted as a force, permeating $\np(d)$ and keeping it from collapsing into a point.

\subsection{Space of normal polytopes}\label{Spaces}
The geometric realiziation of the order simplicial complex of a poset 
$(X,\le)$ will be denoted by $|X|$ and $H_*(-)$ refers to the integral singular homology.

The homology groups of  $|\np(d)|$ can be attacked with the use of methods of \emph{persistent homology} (\cite{Perea} and the many references therein). Because of the obvious free action of $\Aff(\ZZ^d)$  on $|\np(d)|$, we expect that interesting features of the homology $H_*(|\np(d)|)$ can be revealed by studying the finite sub-posets $\np_\epsilon(d)\subset\np(d)$ with $\epsilon$ not too large, i.e., persistent trends quickly become recognizable as $\epsilon\to\infty$. Here $\np_\epsilon(d)$ is the part of $\np(d)$ within the $\epsilon$-neighborhood of $0\in\RR^d$

We conjecture that, at the very core level, $\np(d)$ features a large scale emergent analytic behavior. A strong evidence for this would be a discovery of isolated elements in $\np(d)$, following the strategy, described in Section \ref{Search}. 

Next we introduce another poset, whose geometric realization contains $|\np(d)|$ as a subcomplex, which is topologically more tractable and, conjecturally, provides a handle on the more rigid and rarefied $\np(d)$.

All our cones are pointed, rational, and finitely generated. 

For a cone $C\subset\RR^d$, the additive submonoid $C\cap\ZZ^d\subset\ZZ^d$ has a unique
minimal generating set: the set of indecomposable elements. It is called the \emph{Hilbert basis} of the cone $C$ and denoted by $\Hilb(C)$ \cite[Chapter 2]{kripo}.

When $\Hilb(C)$ lies in an affine subspace $H\subset\RR^d$ with $0\notin H$, the cone $C$ is called \emph{homogeneous}. A finite point configuration in $\RR^d$ is of the form $\Hilb(C)$ for a homogeneous cone $C$ if an only if there exists a normal polytope $P\subset\RR^d$ such that the sets $\Hilb(C)$ and $P\cap\ZZ^d$ are in bijective correspondence by an affine map.

Let $\cones(d)$ be the set of cones $C\subset\RR^d$. It is a partially ordered set as follows: $D< C$ if and only if there exists a finite sequence of cones in $\cones(d)$ of the form
\begin{align*} &D=C_0\subset C_1\subset\ldots\subset C_n=C,\\ &C_i\cap\ZZ^d=\big(C_{i-1}\cap\ZZ^d\big)+\ZZ_{\ge0}x\ \
\text{for some}\ \ x\in C_i\setminus C_{i-1},\quad i=1,\ldots,n. 
\end{align*} 

The subposet of $\cones(d)$, consisting of the nonzero cones in $\big(\RR^{d-1}\times\RR_{>0}\big)\cup\{0\}$, will be denoted by $\cones(d)^+$.  The \emph{homogenization} embedding 
$$
\np(d-1)\to\cones(d)^+,\quad P\to\RR_{\ge0}(P,1), 
$$
is monotonic, inducing an injective continuous map $|\np(d-1)|\to|\cones(d)^+|$. It is known that, for $d\ge4$, the order in $\np(d-1)$ is \emph{weaker} than the order, induced from $\cones(d)$ \cite[Example 2.4]{Cones}. For every natural number $h$, consider the poset $\cones(d)^{(h)}$ of nonzero cones $C\in\cones(d)^+$ with $\Hilb(C)\subset\RR^{d-1}\times[0,h]$, where the order relation is introduced as in $\cones(d)$, but with the additional requirement that the intermediate cones $C_i$ are also in $\cones(d)^{(h)}$. Observe that the embeddings $\cones(d)^{(h)}\hookrightarrow\cones(d)^{(h+1)}$ are not necessarily subposet embeddings. However, being monotonic, they still induce simplicial embeddings of the geometric realizations. 

The following filtration is stable under the natural $\Aff(\ZZ^{d-1})$-action:
$$
\xymatrix{
\np(d-1)=\cones(d)^{(1)}\hookrightarrow\cones(d)^{(2)}\hookrightarrow&\ldots\ar[rr]_{\text{colim}}&&\cones(d)^+.
}
$$

It is easily seen that $\cones(d)^+$ has neither maximal nor nontrivial minimal elements. Informally, $|\cones(d)^+|$ is a result of `smoothing out peaks\rq{} in $|\np(d-1)|$. 

\begin{conjecture}\label{conj-cones}(\cite{Cones}) For every $d,i,j$:
\begin{enumerate}[\rm(a)]
\item
The order relation in $\cones(d)$ is the inclusion order and, thus, $|\cones(d)^+|$ is contractible;
\item
The relative groups
$H_i\big(|\cones(d)^{(j+1)}\setminus\{0\}|,|\cones(d)^{(j)}\setminus\{0\}|\big)$ 
are finitely generated $\ZZ\big[\Aff_{d-1}(\ZZ)\big]$-modules.
\end{enumerate}
\end{conjecture}

Conjecture \ref{conj-cones}(a) for the first nontrivial case $d=3$ is proved in \cite{Cones}.

For $i=0$ and all $j$, Conjecture \ref{conj-cones}(b) in particular implies that the elusive isolated elements in $\np(d)$ form a structured family: for every $j$, after passing from $\cones(d)^{(j)}$ to $\cones(d)^{(j+1)}$, only finitely many such isolated elements (up to unimodular equivalence) cease to be isolated, and all isolated elements are taken out as $j\to\infty$. More generally, the part (b) bounds the complexity of $|\np(d)|$, especially if explicit finite $\ZZ\big[\Aff_{d-1}(\ZZ)\big]$-module generating sets of the relative homology groups can be found.

\subsection{In search of special extremal polytopes}\label{Search} Here we propose three searches, intertwined with deep conjectural properties of $\np(d)$, not even fully explicable at this point.

\medskip\noindent{\bf\emph{Isolated normal polytopes.}} Let  $\Pi=\{(P_{i-1},P_i)\}_\NN$ be a sequence of successive quantum jumps in $\np(d)$. Define the \emph{zeta-value} of $\Pi$ by 
$\zeta_\Pi=\sum_i\vol(P_i)^{-1}$ (possibly $=\infty$).
Here $\vol(-)$ is the normalized volume in $\RR^d$ with respect to $\ZZ^d$. In view of the greedy search for maximal elements in $\np(d)$, discussed in Section \ref{Poset}, the condition $\zeta_\Pi<\infty$ corresponds to an increasingly chaotic growth of normal polytopes via quantum jumps.

Assume Conjecture \ref{conj-cones}(a) is true. Then the sub-poset $\cones(d+1)^{\bf f}\subset\cones(d+1)$ of the full-dimensional cones in $\RR^{d+1}$ is isomorphic to the poset of open subsets $\inte(C)\cap S^d$ of the standard $d$-sphere $S^d$, ordered by inclusion, where $C\subset\RR^{d+1}$ runs over the $(d+1)$-cones and $\inte(-)$ refers to the interior. By \emph{Leray\rq{}s Nerve Theorem} \cite{Nerve}, applied to the open cover of $S^d$ by these open sets, the space $|\cones(d+1)^{\bf f}|$ is homotopic to $S^d$. It is likely that the homotopy equivalence $|\cones(d+1)^{\bf f}|\to S^d$ can be chosen in such a way that the image of the subspace 
$$
|\cones(d+1)^{\bf f}\cap\cones(d+1)^+|\subset|\cones(d+1)^{\bf f}|
$$
is the upper open hemisphere $(S^d)^+\subset S^d$, eventually leading to a homotopy equivalence 
$$
f:|\cones(d+1)^{\bf f}\cap\cones(d+1)^+|\to\RR^d.
$$

Call the sequece $\Pi$ \emph{of Cauchy type} if the images under $f$ of the members of $\Pi$ in $\RR^d$ form a Cauchy sequence in the standard metric.

Now assume $\Pi$ is of Cauchy type and, also, $\zeta_\Pi<\infty$. Assume there exists $Q\in\np(d)$, such that 
$f(Q)=\lim_{n\to\infty}f(P_n)$. (Notice, such a convergence cannot take place on the level of $|\cones(d)|$ because the latter contains $\np(d)$ as a discrete subset.) Then the informal physical interpretation of quantum jumps in Section \ref{Poset} makes $Q$ into a reasonable candidate for an isolated element of $\np(d)$.

An actual computer search can be carried out as follows: 
\begin{enumerate}[\rm$\centerdot$]
\item Find an explicit map $f$ as above; this is essentially achieved in the proof of \cite[Theorem 5.1]{Cones}, where, assuming Conjecture \ref{conj-cones}(a), $\cones(d+1)\setminus\{0\}$ is algorithmically represented as a filtered union of sub-posets $X_i$ ($i\in\NN$) so that each $|X_i|$ contains $S^d$ as a strong deformation retract  (the retractions do not form a compatible system though);
\item
Detect a long sequence $\Pi'=\{(P_{i-1},P_i\}_{i=1}^n$ of successive quantum jumps in $\np(d)$, exhibiting a sharp volume increase together with the Cauchy condition trend for the points $f(P_i)\subset\RR^d$; if one uses the finite approximations to $f$ from \cite[\S5]{Cones} then $\Pi\rq{}$ needs to be chosen correspondingly;
\item
Test elements $Q\in\np(d)$ from $f^{-1}(U)$ on being isolated, where $U\subset\RR^d$ is a small open $d$-ball, centered at $f(P_n)$ and containing $f(P_i)$ for many values of $i$.  
\end{enumerate}

An actual explication of the conjectural large scale emergent analytic behavior of $\np(d)$ may involve the full-blown \emph{polytopal zeta-funtions} $\zeta_\Pi(s)=\sum_{i}\vol(P_i)^{-s}$, 
specializing to the Riemann zeta function in dimension one.

\medskip\noindent{\bf\emph{True random polytopes}.} The successful random search for maximal elements in $\np(d)$, mentioned in Section \ref{Poset}, motivates the approach below.

There are open sources for true random number sequences, e.g., \cite{quant,atmo}. Let  $\Sigma$ be such a sequence, for simplicity assumed to be a $0/1$ sequence. Fix numbers $d,v\in\NN$ and a function $\phi:\NN\to\NN$, and do the following in successive steps. The $n$-th step uses the part $\Sigma(n)$ of $\Sigma$, not used in the previous steps (with the convention $\Sigma(0)=\Sigma$),  as follows: split off the initial segment $\Sigma_n$ from $\Sigma(n)$ of length $n\times d\times v\times \phi(n)$. What is left from $\Sigma(n)$ will be $\Sigma(n+1)$. Partition $\Sigma_n$ in into $d\times v\times\phi(n)$ many $n$-tuples. Each of the $n$-tuples represents a non-negative integer in base $2$. Thus, we have a sequence of length $d\times v\times\phi(n)$ of non-negative integers. Partitioning the latter into $v\times\phi(n)$ many $d$-tuples, we get a sequence  of length $v\times\phi(n)$ of points in $\ZZ_{\ge0}^d$. Partition these sequence of points into $\phi(n)$ many $v$-tuples. For each of these $v$-tuple, compute their convex hull in $\RR^d$. This way we get a sequence of $\phi(n)$ lattice polytopes with at most $v$ vertices. As $n$ grows, we get a sequence of lattice polytopes with at most $v$ vertices. The polytopes are grouped in clusters, with the $n$-th cluster of size $\phi(n)$. What kind of function $\phi(n)$ would be an optimal choice? It cannot be large relative to the number of $0/1$ sequences of length $n$, i.e., $2^n$. Our suggestion is to use $\phi(n)=(ndv)^c$ for some $c\in\NN$; the parameter $c$ can vary. This choice is based on the expectation that, in large scale computations, true randomness is exponentially better than pseudo-randomness.

True random lattice polytopes can be used to determine the frequency of hitting the following objects for $d$ not too large:
\begin{enumerate}[\rm$\centerdot$]
\item Elements of $\np(d)$,
\item Minimal elements of $\np(d)$,
\item Maximal elements of $\np(d)$,
\item Isolated elements of $\np(d)$.
\end{enumerate}

\medskip\noindent{\bf\emph{Random pairs of unimodular simplices.}} Let $\{\ee_1,\ldots,\ee_d\}\subset\ZZ^d$ be the standard basis. The original counterexample to ($\ucp$) \cite{unico} is the convex hull of two parallel unimodular 4-simplicies in $\RR^5$, on lattice distance one apart. In particular, up to $\Aff(\ZZ^5)$-action, it is of the form
\begin{align*}
\conv\big(0,\ee_1,\ldots,\ee_5,(\zz_1,1),\ldots,(\zz_4,1)\big)\subset\RR^5,
\end{align*}
where $\{\zz_1,\ldots,\zz_4\}\subset\ZZ^4$ is a basis. 

Experimentation with randomly generated polytopes of the from
\begin{equation}\label{Hexagon}
\begin{aligned}
\conv\big(0,\ee_1,\ldots,\ee_{d+1},(\zz_1,1),\ldots,(\zz_d,1)\big)\subset\RR^{d+1},
\end{aligned}
\end{equation}
where $\{\zz_1,\ldots,\zz_d\}\subset\ZZ^d$ is a basis, did not yield new counterexamples to ($\icp$). We suggest to enhance this search by incorporating true randomness, combined with a fast basis generation in $\ZZ^d$, using \cite{Carter}. In more detail, it is shown in \cite{Carter} that every matrix $B\in SL_d(\ZZ)$ factors into $\theta(d)=36+\frac12(3d^2-d)$ elementary matrices $e_{ij}^a$: 1-s on the main diagonal and the only off-diagonal non-zero entry $a$ in the position $ij$. In particular, for every polytope of the form (\ref{Hexagon}), we can assume that $\zz_1,\ldots,\zz_d$ are the rows of a matrix
\begin{align*}
\prod_{k=1}^{\theta(d)} e_{i_kj_k}^{a_k},\quad1\le i_k\not=j_k\le d,\quad a_k\in\ZZ.
\end{align*}
We get a parametrization of the polytopes of type (\ref{Hexagon}) by the sequences
\begin{align*}
(a_1,\ldots,a_{\theta(d)},i_1,j_1,\ldots,i_{\theta(d)},j_{\theta(d)}),&\\
\qquad\qquad\qquad1\le i_k\not=j_k\le d,&\\
a_k\in\ZZ&\qquad(k=1,\ldots,\theta(d)).
\end{align*}

For the first $\theta(d)$ components one can use true random numbers and, for not too large values of $d$, the last $2\theta(d)$ components can be generated systematically.

The generated polytopes of type (\ref{Hexagon}) can be checked on being normal, minimal, maximal, violating ($\ucp$) -- in this succession. The most expensive Carath\'eodory rank computations only makes sense for the detected counterexamples to ($\ucp$). 

\section{Back to continuous polytopes}\label{Pyramidal} 

The poset of general convex polytopes in $\RR^d$, ordered by inclusion, is topologically trivial. However, its discrete model $\np(d)$ and the conic enhancement $\cones(d+1)$ of the latter suggest other natural orders on the set of convex polytopes and the question whether these orders coincide with the inclusion order is a challenge in polytope theory of classical flavor. 

The set of polytopes in $\bigoplus_\NN\RR$ of dimension at most $d$ will be denoted by $\POL(d)$. The set of all polytopes will be denoted by $\POL(\infty)$. For a subfield $\k\subset\RR$, the corresponding sets of polytopes with vertices in $\oplus_\NN\k$ will be denoted by $\POL_\k(d)$ and $\POL_\k(\infty)$. 

The \emph{Hausdrorff distance} between two nonempty compact subsets $X,Y\subset\oplus_\NN\RR$ will be denoted by $\d(X,Y)$ \cite[Chapter 1.2]{Grunbaum}. For a sequence of polytopes $\{P_i\}_\NN$ and a polytope $Q$, we write $\underset{i\to\infty}\lim P_i=Q$ if $\underset{i\to\infty}\lim\d(P_i,Q)=0$.

Let $P$ be a polytope. A \emph{pyramid over} or \emph{with base} $P$ is the convex hull $Q$ of $P$ and a point $v$, not in the affine hull of $P$. The point $v$ is the \emph{apex} of $Q$.

Let $\k\subset\RR$ be a subfield, $d\le\infty$, and $\POL$ be one of the sets $\POL(d)$, $\POL_{\k}(d)$.

\begin{definition}\label{Pyrmidal_growth}
(a) A pair of polytopes $P\subset Q$ in $\POL$ forms a \emph{pyramidal extension} if $\Delta=\overline{Q\setminus P}$ is a pyramid and $\Delta\cap P$ is a facet of $\Delta$, i.e., either $Q$ is a pyramid over $P$ or obtained from $P$ by stacking a pyramid onto a facet. For a pyramidal extension $P\subset Q$ we write $P\D\subset Q$. (The bar `$^{\overline{\ \ \ }}$\rq{} refers to the topological closure.)

\noindent(b) The partial order on $\POL$, generated by the pyramidal extensions within $\POL$, will be denoted by $\D\le$ and called the \emph{pyramidal growth}.

\noindent(c) The \emph{transfinite pyramidal growth} $\Infty\le$ is the smallest partial order  on $\POL$, containing $\D\le$ and satisfying $P\Infty\le Q$ whenever there exists an ascending sequence  $P=P_0\Infty\le P_1\Infty\le P_2\Infty\le\ldots$ with  $Q=\underset{i\to\infty}\lim P_i$. 
\end{definition}

Pyramidal extensions are more general than the extensions, used in the definition of \emph{stacked polytopes} \cite[Chapter 10.6]{Grunbaum}, along with their direct generalization to arbitrary initial polytopes: by only allowing the stackings of pyramids onto facets when none of the codimension 2 faces disappears we get a different partial order on polytopes. Obviously, it does not coincide with the inclusion order. We do not know whether the transfinite completion of this order is the same as $\Infty\le$.

The motivation for introducing the pyramidal growth, apart from being a natural geometric construction, is twofold:

\begin{enumerate}[\rm$\centerdot$]
\item It is a continuous counterpart of the special/minimal extensions $C\subset D$ in $\cones(d)$, when the new extremal generator $v\in D\setminus C$ is on lattice distance 1 from the affine hull of every facet $F\subset C$ that separates $v$ and $\inte(C)$. We do not know whether these types of cone extensions already suffice to generate the inclusion order for cones. The positive answer to this question would prove Conjecture \ref{conj-cones} and it was used to develop a computational evidence for that conjecture (\cite{Chan} and the reference at the end of \cite[\S D2]{Cones}). 
\item The poset $\big(\POL_\QQ(\infty),\D\le\big)$ is implicit in our $K$-theoretic works on monoid rings. Informally, a pyramidal extension $P\D\subset Q$ represents a minimal enlargement of a polytope, allowing to transfer certain $K$-theoretic information on $P$ to the polytope $Q$. More precisely, rational polytopes give rise to  submonoids of $\ZZ^d$ and, when $P\D\subset Q$, certain $K$-theoretic objects over the monoid ring, associated with $Q$, are extended from the submonoid ring, associated with $P$. Results of this type in various $K$-theoretic scenarios are obtained in \cite{GuAND,Nilpotence,Elrows}.
\end{enumerate}

\medskip A continuous counterpart of the order relation in $\cones(d)$ is given by the following

\begin{definition}\label{Quasi_growth}
(a) A \emph{quasi-pyramidal growth} in $\POL$, is a pair of polytopes $P\subset Q$, admitting within $\POL$ a finite sequence of polytopes
\begin{align*}
P=P_0\subset P_1\subset\ldots\subset P_n=Q
\end{align*}
and pyramidal extensions
$$
P'_i\D\subset P_i,\qquad i=1,\ldots,n,
$$
such that 
$$
P_0\subset P_i'\subset P_{i-1},\qquad i=1,\ldots,n.
$$
The resulting partial order on $\POL$ will be denoted by $\q\le$. 

\noindent(b) The \emph{quasi-pyramidal defect} of a pair of polytopes $P\q\le Q$ in $\POL$ is defined by
\begin{align*}
\delta(P,Q)=\inf\bigg\{\sum \d(P_i',P_i)\ \bigg |\ P'_i\subsetneq P_{i-1}\bigg\}, 
\end{align*}
where the infimum is taken over the sequences as in the part (a).
\end{definition}

\begin{theorem}\label{Pyramidal}(\cite{Pyramidal})
Let $\k\subset\RR$ be a subfield.
\begin{enumerate}[\rm(a)]
\item $\Infty\le$ is the inclusion order on $\POL_\k(\infty)$.
\item $\D\le$ is the inclusion order on $\POL_\k(3)$.
\item Assume $\D\le$ is the inclusion order on $\POL_\k(d)$. Then $\q\le$ is the inclusion order on $\POL_\k(d+1)$ and $\delta(P,Q)=0$ for any two polytopes $P\subset Q$ in $\POL_\k(d+1)$.
\end{enumerate}
\end{theorem}

In particular, $\q\le$ is the inclusion order on $\POL(4)$.

We expect that Theorem \ref{Pyramidal}(b) is sharp and Theorem \ref{Pyramidal}(c) holds unconditionally in all dimensions, in parallel with Conjecture \ref{conj-cones}(a):

\begin{conjecture}\label{continuous}
 \begin{enumerate}[\rm(a)] 
 \item The pyramidal growth is not the inclusion order in $\POL(4)$;
 \item For every $d$, the quasi-pyramidal growth is the inclusion order in $\POL(d)$ and $\delta(P,Q)=0$ for any polytopes $P\subset Q$.
 \end{enumerate}
  \end{conjecture} 

We propose two possible strategies for proving the part (a) of this conjecture. 

If $P\D\subset Q$ and $P$ is a rational polytope, then $Q$ is combinatorially equivalent to a rational polytope. Whether a relation $P\D\le Q$ imposes a constrain on the field over which $Q$ is defined, is an interesting question. In view of the existence of 4-dimensional combinatorial types, whose geometric realizations require any pre-specified real finite field extension $\QQ\subset\k$ \cite{Realization}, the positive answer to this question would imply that $\D\le$ and $\subset$ are different. 

The other strategy is related to \emph{approximations by polytopes} \cite{Gruber}. It is not clear whether, for $d\ge4$, there exists a sequence $0=P_0\ \D\le\ P_1\ \D\le\ \ldots$ in $\POL(d)$, starting with the origin $0\in\RR^d$ and such that $\bigcup_{i=0}^\infty P_i$ is the open unit $d$-ball. If this is not the case then $\D\le$ is not the inclusion order on $\POL_\QQ(d)$.

\bigskip\noindent\emph{Acknowledgment.} I thank Winfried Bruns for the many years of our collaboration, which extends beyond the topics discussed here. Apart of the mathematical collaboration, he implemented general algorithms, resulting from this collaboration, and  carried out decisive computer experimentations. For the all too speculative ideas above I take the sole responsibility.  
\bibliography{bibliography}
\bibliographystyle{plain}

\end{document}